\documentclass[12pt]{article}

\usepackage[T1]{fontenc}
\usepackage{amsmath}
\usepackage{amssymb}
\usepackage{diagrams}
\usepackage[noBBpl]{mathpazo}

\renewcommand{\texttt}[1]{{\fontfamily{pcr}\fontseries{m}\fontshape{n}\selectfont#1}}

\providecommand{\qed}{\hspace*{\fill}$\Box$}
\newcommand{\isopil}{\stackrel{\raisebox{0.1ex}[0ex][0ex]{\(\sim\)}}%
			{\raisebox{-0.15ex}[0.28ex]{\(\rightarrow\)}}}
\newcommand{\tensor}	{\otimes}
\newcommand{\df}{\: {\raisebox{0.255ex}{\normalfont\scriptsize :\!\!}}=}

\providecommand{\overskrift}[1]{\par\noindent\relax{\LARGE #1}\par\bigskip}

\newcommand{\hovedfont}{\normalfont\bfseries}

\usepackage{theorem}
\theorembodyfont{\normalfont\itshape}
\theoremheaderfont{\hovedfont}
	\theoremstyle{change}
\newtheorem{lemma}{Lemma.}[section]
\newtheorem{prop}[lemma]{Proposition.}

\newtheorem{cor}[lemma]{Corollary.}
	\theorembodyfont{\normalfont}

\newtheorem{taller}[lemma]{$\!\!$}
\newenvironment{blanko}[1]%
{\begin{taller}{\hovedfont #1}\normalfont}%
{\end{taller}}

\newenvironment{dem}%
{%
\begin{list}{\em Proof. }%
{\setlength{\labelsep}{0mm}\setlength{\leftmargin}{0mm}%
\setlength{\labelwidth}{0mm}\setlength{\listparindent}{\parindent}%
\setlength{\parsep}{\parskip}\setlength{\partopsep}{0mm}}%
\item%
}%
{%
\qed\end{list}%
}
\newenvironment{dem*}[1]%
{%
\begin{list}{\em #1 }%
{\setlength{\labelsep}{0mm}\setlength{\leftmargin}{0mm}%
\setlength{\labelwidth}{0mm}\setlength{\listparindent}{\parindent}%
\setlength{\parsep}{\parskip}\setlength{\partopsep}{0mm}}%
\item%
}%
{%
\qed\end{list}%
}

\newenvironment{blanko*}[1]%
{%
\begin{list}{\bf {#1} }%
{\setlength{\labelsep}{0mm}\setlength{\leftmargin}{0mm}%
\setlength{\labelwidth}{0mm}\setlength{\listparindent}{\parindent}%
\setlength{\parsep}{\parskip}\setlength{\partopsep}{0mm}}%
\item%
}%
{%
\end{list}%
}

\usepackage{mathrsfs}
\newcommand{\CC}{\mathscr{C}}

\providecommand{\kat}[1]{\text{\textbf{\textsl{#1}}}}
\newcommand{\id}{\operatorname{id}}
\newcommand{\Cat}{\kat{Cat}}

\makeatletter
\renewcommand{\ps@headings}
  {\setlength{\headheight}{12pt}%
   \setlength{\headsep}{11pt}%
   \renewcommand{\@oddhead}{\parbox{\textwidth}{%
      \scriptsize
      \texttt{file: commutativity.tex \ \ \ \ version: 2006-08-15 22:52
      \hfill page: [\thepage/\pageref{lastpage}]}
      \\ \rule[8pt]{\textwidth}{0.3pt}}%
   }
  \renewcommand{\@oddfoot}{}
  \renewcommand{\@evenfoot}{}%
}
\makeatother

\newcommand{\halfgrid}{%
\put(0,0){\line(1,0){40}}%
\put(0,20){\line(1,0){40}}%
\put(0,40){\line(1,0){40}}%
\put(0,0){\line(0,1){40}}%
\put(20,0){\line(0,1){40}}%
\put(40,0){\line(0,1){40}}%
}

\newcommand{\twobytwogrid}{%
  \put(0,0){\line(1,0){20}}
  \put(0,10){\line(1,0){20}}
  \put(0,20){\line(1,0){20}}
  \put(0,0){\line(0,1){20}}
  \put(10,0){\line(0,1){20}}
  \put(20,0){\line(0,1){20}}
}

\begin{document}  

\def\vspec#1{\special{ps:#1}}
\def\rotstart#1{\vspec{gsave currentpoint currentpoint translate
	#1 neg exch neg exch translate}}
\def\rotfinish{\vspec{currentpoint grestore moveto}}
	
\def\psrotate#1#2{\rotstart{#1 rotate}#2\rotfinish}

\def\psr#1{\rotstart{36.87 rotate}\hbox to0pt {\vsize 0pt \hss\(\scriptstyle #1\)\hss}\rotfinish}

\def\PSR#1#2{\rotstart{#1 rotate}\hbox to0pt {\vsize 0pt \hss\(#2\)\hss}\rotfinish}

\newcommand{\verticalSimeq}{\PSR{90}{\simeq}}

\pagestyle{headings}

\vspace*{10pt}
\begin{center}
  
  \overskrift{Note on commutativity in double semigroups \\[6pt]
  and two-fold monoidal categories}

\medskip

\noindent
\textsc{Joachim Kock}
\smallskip

\normalsize

\end{center}

\begin{abstract}
  A concrete computation --- twelve slidings with sixteen tiles ---
  reveals that certain commutativity phenomena occur in every double
  semi\-group.  This can be seen as a sort of Eckmann-Hilton argument,
  but it does not use units.  The result implies in particular that
  all cancellative double semi\-groups and all inverse double semi\-groups
  are commutative.  Stepping up one dimension, the result is used to
  prove that all strictly associative two-fold monoidal categories 
  (with weak units) are
  degenerate symmetric.  In particular, strictly associative
  one-object, one-arrow \linebreak $3$-groupoids (with weak units)
  cannot realise all simply-connected homotopy \linebreak $3$-types.
\end{abstract}


\section{Introduction and results}

\begin{blanko}{The Eckmann-Hilton argument.}
  In 1932, \v{C}ech introduced the higher homotopy groups $\pi_i$,
  $i>1$ in a contribution submitted to the International Congress of
  Mathematicians in Z\"urich.  His paper was received by Alexandrov
  and Hopf who quickly realised that all these groups are abelian (or 
  perhaps \v Cech had noticed this himself), and
  for this reason they felt it could not be the correct notion.  They
  convinced \v{C}ech to withdraw his paper, and in the final
  proceedings only a very short communication of \v{C}ech was 
  included \cite{Cech:1932}.%
  \footnote{I learned this story from Ronnie Brown, who in turn got it
  from Eldon Dyer.}

  The natural generality of the commutativity argument, known as the
  Eckmann-Hilton argument \cite{Eckmann-Hilton}, is an elementary
  statement about double monoids (although asssociativity is not 
  essential).  Recall that a double monoid is a
  set $S$ equipped with two compatible monoid structures, i.e.~two
  associative and unitary multiplications $*_h$ and $*_v$ satisfying
  the interchange law
  $$
  (x *_h y) *_v (z*_h w) = (x *_v z) *_h (y *_v w)
  $$
  for all $x,y,z,w \in S$.
  If $*_h$ is depicted
  horizontally and $*_v$ vertically,
%
    then the interchange law says that in a composite composition
  \begin{center}
  \setlength{\unitlength}{1.5pt}
  \begin{picture}(20,20)(0,0)
     \twobytwogrid
    \put(5,3.5){\makebox(0,0)[b]{\scriptsize $x$}}
    \put(15,2){\makebox(0,0)[b]{\scriptsize $y$}}
    \put(5,13){\makebox(0,0)[b]{\scriptsize $z$}}
    \put(15,13){\makebox(0,0)[b]{\scriptsize $w$}}
  \end{picture}
\end{center} 
  it doesn't matter whether the vertical or the horizontal composition 
  is performed first.  
  
    It is automatic from the interchange law that the two units
  coincide.  This unit, which we denote by $1$, obviously commutes
  with every other element, for the horizontal as well as the vertical
  composition.  This innocent-looking special case of commutativity in
  fact forces the two composition laws to be commutative, and to
  coincide.  This is the Eckmann-Hilton argument:
  \begin{center}
  \setlength{\unitlength}{1.5pt}
  \begin{picture}(20,10)(0,-5)
    \put(0,0){\line(0,1){10}}
    \put(10,0){\line(0,1){10}}
    \put(20,0){\line(0,1){10}}
    \put(0,0){\line(1,0){20}}
    \put(0,10){\line(1,0){20}}
    \put(5,3){\makebox(0,0)[b]{\scriptsize $a$}}
    \put(15,3){\makebox(0,0)[b]{\scriptsize $b$}}
  \end{picture}
  \begin{picture}(24,20)(0,0)
    \put(12,10){\makebox(0,0){$=$}}
  \end{picture}
  \begin{picture}(20,20)(0,0)
     \twobytwogrid
    \put(5,3){\makebox(0,0)[b]{\scriptsize $a$}}
    \put(15,13){\makebox(0,0)[b]{\scriptsize $b$}}
    \put(5,13){\makebox(0,0)[b]{\scriptsize $1$}}
    \put(15,3){\makebox(0,0)[b]{\scriptsize $1$}}
  \end{picture}
  \begin{picture}(24,20)(0,0)
    \put(12,10){\makebox(0,0){$=$}}
  \end{picture}  
  \begin{picture}(10,20)(0,0)
    \put(0,0){\line(1,0){10}}
    \put(0,10){\line(1,0){10}}
    \put(0,20){\line(1,0){10}}
    \put(0,0){\line(0,1){20}}
    \put(10,0){\line(0,1){20}}
    \put(5,3){\makebox(0,0)[b]{\scriptsize $a$}}
    \put(5,13){\makebox(0,0)[b]{\scriptsize $b$}}
  \end{picture}
  \begin{picture}(24,20)(0,0)
    \put(12,10){\makebox(0,0){$=$}}
  \end{picture}
  \begin{picture}(20,20)(0,0)
     \twobytwogrid
    \put(15,3){\makebox(0,0)[b]{\scriptsize $a$}}
    \put(5,13){\makebox(0,0)[b]{\scriptsize $b$}}
    \put(5,3){\makebox(0,0)[b]{\scriptsize $1$}}
    \put(15,13){\makebox(0,0)[b]{\scriptsize $1$}}
  \end{picture}
  \begin{picture}(24,20)(0,0)
    \put(12,10){\makebox(0,0){$=$}}
  \end{picture}
  \begin{picture}(20,10)(0,-5)
    \put(0,0){\line(0,1){10}}
    \put(10,0){\line(0,1){10}}
    \put(20,0){\line(0,1){10}}
    \put(0,0){\line(1,0){20}}
    \put(0,10){\line(1,0){20}}
    \put(5,3){\makebox(0,0)[b]{\scriptsize $b$}}
    \put(15,3){\makebox(0,0)[b]{\scriptsize $a$}}
  \end{picture}
\end{center} 

%

  (Since the homotopy group $\pi_2$ is defined in terms of maps from
  squares, which can be composed in two compatible ways (horizontally
  and vertically), the commutativity of $\pi_2$ as well as all higher
  homotopy groups follows readily \cite{Eckmann-Hilton}.)
\end{blanko}

\begin{blanko}{Double semigroups.}
  It is clear that the existence of the unit is a key ingredient in
  the Eckmann-Hilton argument.  For double semigroups (i.e.~sets with
  two compatible non-unital associative multiplications), the argument
  does not work, and indeed it is easy to exhibit examples of double
  semigroups which are not commutative.
  
  The main discovery of the present note is that certain commutativity
  phenonema do arise even in double semigroups without units.  These
  phenomena occur in expressions with many terms, where the interchange
  law and associativity combined give rise to some rearrangements of
  terms.  One such commutativity is expressed by
  
\begin{quote}
\begin{blanko*}{Proposition~\ref{ab}.}{\em
  For any sixteen elements $a,b,\ldots$ in any double semigroup, this
  equation holds:

  \begin{center}
  \setlength{\unitlength}{1.5pt}
  \begin{picture}(40,40)(0,0)
    \halfgrid
    \put(0,10){\line(1,0){40}}
    \put(0,30){\line(1,0){40}}
    \put(10,0){\line(0,1){40}}
    \put(30,0){\line(0,1){40}}
    \put(15,23){\makebox(0,0)[b]{\scriptsize $a$}}
    \put(25,23){\makebox(0,0)[b]{\scriptsize $b$}}
  \end{picture}
  \begin{picture}(20,40)(0,0)
    \put(10,20){\makebox(0,0){$=$}}
  \end{picture}
  \begin{picture}(40,40)(0,0)	
    \halfgrid
    \put(0,10){\line(1,0){40}}
    \put(0,30){\line(1,0){40}}
    \put(10,0){\line(0,1){40}}
    \put(30,0){\line(0,1){40}}
    \put(15,23){\makebox(0,0)[b]{\scriptsize $b$}}
    \put(25,23){\makebox(0,0)[b]{\scriptsize $a$}}
  \end{picture}
\end{center}
(The empty boxes represent fourteen nameless elements, the same on
each side of the equation, and in the same order.)}
\end{blanko*}
\end{quote}

  The proof is an elementary computation exploiting the geometrical
  representation of the two multiplication laws.
  It would be quite cumbersome to write it out algebraically.
\end{blanko}
  
\begin{blanko}{Cancellative double semigroups and inverse double
  semigroups are commutative.} In double semigroups with some further
  cancellation properties, the blanks in the equation can be
  cancelled away and we get absolute commutativity.  In particular,
  every cancellative double semigroup is commutative 
  (Corollary~\ref{cancellative}).
  It is also shown that every inverse double semigroup is
  commutative (Proposition~\ref{inverse}). 
\end{blanko}

\begin{blanko}{Two-fold monoidal categories and braidings.}
  The notion of double monoid as well as the Eckmann-Hilton argument
  make sense in any monoidal category in place of the category of
  sets.  A double monoid in $\Cat$ is the same thing as a category
  with two compatible strict monoidal structures, and the
  Eckmann-Hilton argument shows that such are commutative.  
  
  It is natural to consider the non-strict version of this situation,
  i.e.~a category equipped with two (non-strict) monoidal structures
  which are compatible up to coherent isomorphism.  In this case the
  Eckmann-Hilton argument consists of a sequence of specific
  isomorphisms, and it turns out to define a braiding.  This was
  observed by Joyal and Street~\cite{Joyal-Street:braided-tensor} in
  1985, and in fact was one of the motivating examples for their
  discovery of the notion of braided monoidal category.  (Conversely a
  braiding on a monoidal category can be used to construct a second,
  weakly compatible, monoidal structure.)
  
  In this way, the Eckmann-Hilton argument is
  directly related to core subtleties of higher category theory.
\end{blanko}

\begin{blanko}{Strictifications.}
  One of the key themes in higher category theory is strictification.
  Finding strict or semi-strict models for weak structure often
  amounts to powerful coherence results.  It is well known that every
  monoidal category is equivalent to a strict one
  \cite{MacLane:naturalAssociativity} \cite[XI.3]{MacLane:categories},
  but the argument of the previous paragraph shows that not every
  two-fold monoidal category is equivalent to a strict one --- which
  is just another expression of the fact that not every braided
  monoidal category is braided equivalent to a commutative monoid in
  $\Cat$.  In fact, two-fold monoidal categories can be seen as
  tri\-categories with one object and one arrow, and the observation is
  the simplest case of the fact that not all tricategories are
  equivalent to strict $3$-categories~\cite{Gordon-Power-Street}.
  
  As a rule of thumb, one can strictify one level of weak structure,
  but in general not two levels at the same time; see for example
  Paoli~\cite{Paoli:0607} who studies two different one-level
  strictifications of weak $3$-groupoids in the sense of
  Tamsamani~\cite{Tamsamani:thesis}, and compares with
  $\operatorname{cat}^2$-groups \cite{Loday:catn} in the path-connected
  case.  However, it is sometimes possible to strictify one level of
  structure and parts of other levels.  For example, every two-fold 
  monoidal category is equivalent to a two-fold strict
  monoidal category with iso-interchange~\cite{Balteanu-et.al.:9808} 
   (corresponding to the fact that every braided monoidal
   category is equivalent to a braided strict monoidal category).
  
  The braiding that results from the $\Cat$ version of the
  Eckmann-Hilton argument is a composite of unit structure
  isomorphisms and interchange isomorphisms.  As just explained, it is
  possible to strictify the units if the interchange law is kept weak.
  It is a natural question whether it would be possible to strictify
  instead the interchange law while keeping the unit weak.  This idea
  is related to Simpson's conjecture:
\end{blanko}  

\begin{blanko}{Simpson's conjecture.}
  Based on a careful analysis of strict $3$-groupoids, and the
  observation that the units play a key role in the Eckmann-Hilton
  argument, Simpson~\cite{Simpson:9810} was led to suspect that units
  ( = identity arrows)
  can account for all higher homotopical data in higher categories.  A
  strong version of his conjecture states roughly that every weak
  $n$-category is equivalent to one where only the units are weak.
  (See \cite{Simpson:9810} and \cite{Kock:0507} for more formal
  statements of the conjecture.)  Simpson's conjecture is highly
  surprising and goes against all trends in higher category theory,
  where the emphasis was always on the composition laws.
  
  A weaker version of the conjecture states that strictly associative
  $n$-groupoids with weak units is a model for homotopy $n$-types.
  In contrast, completely strict $n$-groupoids can model only 
  homotopy $n$-types whose higher Whitehead brackets (i.e., those 
  beyond the action of $\pi_1$) are zero.
  The first interesting case is that of $1$-connected homotopy
  $3$-types, since this is the first appearance of a non-trivial
  higher Whitehead bracket,\break $\pi_2 \tensor \pi_2 \to \pi_3$,
  which in turn is closely related to the braiding.  In fact,
  $1$-connected $3$-types can be realised by braided categorical
  groups, cf.~Joyal-Tierney~\cite{Joyal-Tierney:homotopy-types}; see
  also Brown-Gilbert~\cite{Brown-Gilbert} for a closely related model.
  A version of Simpson's conjecture in this
  case was proved in Joyal-Kock~\cite{Joyal-Kock:traintracks}:
  one-object $3$-groupoids which are strict in all respects except
  that there are only weak identity arrows (forming a contractible
  space) can model all $1$-connected homotopy $3$-types.
\end{blanko}
    
\begin{blanko}{Strictification of composition and interchange?}
  It is tempting to reformulate (and distort) Simpson's conjecture by saying that
  composition and interchange can always be strictified if just the
  units remain weak.  One might think that every two-fold monoidal
  category is equivalent to one with strict compositions and strict
  interchange.
%
  However, this is false: in Section~\ref{Sec:cat} of this note, the
  commutativity in Proposition~\ref{ab} is used to show that
  
\begin{quote}
\begin{blanko*}{Proposition~\ref{deg-sym}.}{\em
  If $\CC$ is a two-fold monoidal category,  strictly associative 
  and with strict interchange,
  then $\CC$ is degenerate symmetric (i.e.~has a symmetry $\sigma$ such
that $\sigma_{X,X} = \id_{X\tensor X}$).}
\end{blanko*}
  \end{quote}

  Recall from \cite{Balteanu-et.al.:9808} that braided monoidal
  categories correspond to two-fold loop spaces (via group completion
  of the nerve) and arbitrary $1$-connected homotopy $3$-types
  \cite{Joyal-Street:braided-tensor}, while symmetric monoidal
  categories correspond to infinite loop spaces.  Since infinite loop
  spaces have vanishing Whitehead bracket $\pi_2\tensor\pi_2 \to
  \pi_3$, we find the following corollary.
  
\begin{quote}
\begin{blanko*}{Corollary~\ref{3-groupoids}.}
  {\em Strictly associative one-object, one-arrow $3$-groupoids (but
  still with weak units) cannot realise all simply-connected homotopy
  $3$-types.}
  \end{blanko*}
\end{quote}

%
\end{blanko}

\begin{blanko*}{Acknowledgements.}
  The sliding argument of Proposition~\ref{ab} was discovered after
  conversations with Andr\'e Henriques, whom I thank for precious
  input.  I am also thankful to Ronnie Brown for some pertinent comments.
\end{blanko*}

\section{Double semigroups}
\label{Sec:set}

A {\em double semigroup} is a set equipped with two compatible associative
multiplication laws.  In other words, it is just like a double monoid,
except that there is no unit.

\begin{blanko}{Example.}
  A double semigroup is not necessarily commutative: take any set with at
  least two elements, and let both composition laws be the
  `K-combinator': $x *_v y \df x *_h y \df x $.  Clearly, this is
  associative, and satisfies the interchange law, but it is not commutative.
\end{blanko}

\begin{blanko}{Sliding tiles.}
When writing  the graphical   representation of some product in a 
double semi\-group, it is important to note that there is
a certain freedom in where to set the `walls' --- this comes about 
because of associativity.  For example, in the
product
$$
(c *_h d *_h e) *_v(a *_h b)
$$
there are two ways of setting parentheses in $(c *_h d *_h e)$:
$$
((c *_h d) *_h e) *_v (a *_h b)  = (c *_h (d *_h e)) *_v (a *_h b) .
$$
Hence graphically we get
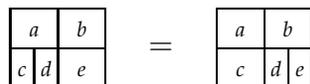
\begin{center}
\setlength{\unitlength}{1.5pt}
\begin{picture}(24,20)(0,0)
  \put(0,0){\line(1,0){24}}
  \put(0,10){\line(1,0){24}}
  \put(0,20){\line(1,0){24}}
  \put(0,0){\line(0,1){20}}
  \put(12,0){\line(0,1){20}}
  \put(24,0){\line(0,1){20}}
  
  \put(6,0){\line(0,1){10}}

  \put(6,13){\makebox(0,0)[b]{\scriptsize $a$}}
  \put(18,13){\makebox(0,0)[b]{\scriptsize $b$}}
  \put(3,3){\makebox(0,0)[b]{\scriptsize $c$}}
  \put(9,3){\makebox(0,0)[b]{\scriptsize $d$}}
  \put(18,3){\makebox(0,0)[b]{\scriptsize $e$}}
\end{picture}
\begin{picture}(24,20)(0,0)
  \put(12,10){\makebox(0,0){$=$}}
\end{picture}
\begin{picture}(24,20)(0,0)
  \put(0,0){\line(1,0){24}}
  \put(0,10){\line(1,0){24}}
  \put(0,20){\line(1,0){24}}
  \put(0,0){\line(0,1){20}}
  \put(12,0){\line(0,1){20}}
  \put(24,0){\line(0,1){20}}
  
  \put(18,0){\line(0,1){10}}

  \put(6,13){\makebox(0,0)[b]{\scriptsize $a$}}
  \put(18,13){\makebox(0,0)[b]{\scriptsize $b$}}
  \put(6,3){\makebox(0,0)[b]{\scriptsize $c$}}
  \put(15,3){\makebox(0,0)[b]{\scriptsize $d$}}
  \put(21,3){\makebox(0,0)[b]{\scriptsize $e$}}
\end{picture}
\end{center}
The upshot is that in the graphical representation, sliding the inner
walls of a given rectangle does not change the corresponding algebraic
expression.  It is clear that such slidings can never change the order
of the elements that touch the border of the expression.  Hence in the
picture above, the order of the five elements will always be
$c,d,e,b,a$, walking around counter-clockwise.  But in expressions
with more elements, nontrivial permutations can take place, as the
following computation shows.
\end{blanko}

\begin{prop}\label{ab}
  For any sixteen elements $a,b,\ldots$ in any double semigroup, this
  equation holds:

  \begin{center}
  \setlength{\unitlength}{1.5pt}
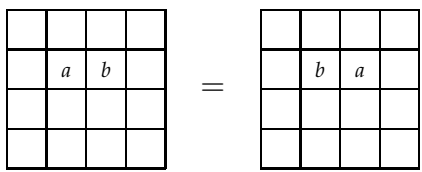
  \begin{picture}(40,40)(0,0)
    \halfgrid
    \put(0,10){\line(1,0){40}}
    \put(0,30){\line(1,0){40}}
    \put(10,0){\line(0,1){40}}
    \put(30,0){\line(0,1){40}}
    \put(15,23){\makebox(0,0)[b]{\scriptsize $a$}}
    \put(25,23){\makebox(0,0)[b]{\scriptsize $b$}}
  \end{picture}
  \begin{picture}(20,40)(0,0)
    \put(10,20){\makebox(0,0){$=$}}
  \end{picture}
  \begin{picture}(40,40)(0,0)	
    \halfgrid
    \put(0,10){\line(1,0){40}}
    \put(0,30){\line(1,0){40}}
    \put(10,0){\line(0,1){40}}
    \put(30,0){\line(0,1){40}}
    \put(15,23){\makebox(0,0)[b]{\scriptsize $b$}}
    \put(25,23){\makebox(0,0)[b]{\scriptsize $a$}}
  \end{picture}
\end{center}
(The empty boxes represent fourteen nameless elements, the same on
each side of the equation, and in the same order.)
\end{prop}

\begin{dem}
  We shall perform twelve slidings, each representing a strict 
  equality.  We only label the middle four elements, since anyway the 
  configuration of the
  elements touching the border is fixed by any sliding.

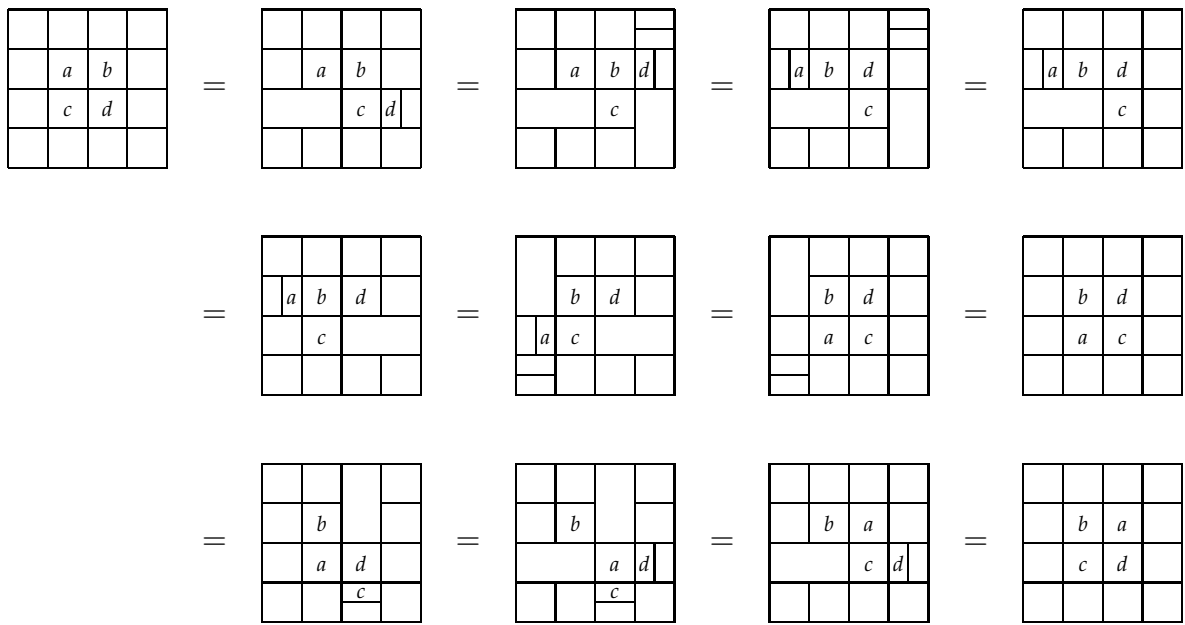
\begin{center}
  \setlength{\unitlength}{1.5pt}
  \begin{picture}(40,40)(0,0)	
    \halfgrid
    \put(0,10){\line(1,0){40}}
    \put(0,30){\line(1,0){40}}
    \put(10,0){\line(0,1){40}}
    \put(30,0){\line(0,1){40}}
    \put(15,23){\makebox(0,0)[b]{\scriptsize $a$}}
    \put(25,23){\makebox(0,0)[b]{\scriptsize $b$}}
    \put(15,13){\makebox(0,0)[b]{\scriptsize $c$}}
    \put(25,13){\makebox(0,0)[b]{\scriptsize $d$}}
  \end{picture}
  \begin{picture}(20,50)(0,0)
    \put(10,20){\makebox(0,0){$=$}}
  \end{picture}
  \begin{picture}(40,40)(0,0)	
    \halfgrid
    \put(0,10){\line(1,0){40}}
    \put(0,30){\line(1,0){40}}
    \put(10,0){\line(0,1){10}}
    \put(10,20){\line(0,1){20}}
    \put(30,0){\line(0,1){40}}
    \put(35,10){\line(0,1){10}}

    \put(15,23){\makebox(0,0)[b]{\scriptsize $a$}}
    \put(25,23){\makebox(0,0)[b]{\scriptsize $b$}}
    \put(25,13){\makebox(0,0)[b]{\scriptsize $c$}}
    \put(32.5,13){\makebox(0,0)[b]{\scriptsize $d$}}
  \end{picture}
  \begin{picture}(20,40)(0,0)
    \put(10,20){\makebox(0,0){$=$}}
  \end{picture}
  \begin{picture}(40,40)(0,0)	
    \halfgrid
    \put(0,10){\line(1,0){30}}
    \put(0,30){\line(1,0){40}}
    \put(30,35){\line(1,0){10}}
    \put(10,0){\line(0,1){10}}
    \put(10,20){\line(0,1){20}}
    \put(30,0){\line(0,1){40}}
    \put(35,20){\line(0,1){10}}

    \put(15,23){\makebox(0,0)[b]{\scriptsize $a$}}
    \put(25,23){\makebox(0,0)[b]{\scriptsize $b$}}
    \put(25,13){\makebox(0,0)[b]{\scriptsize $c$}}
    \put(32.5,23){\makebox(0,0)[b]{\scriptsize $d$}}
  \end{picture}
  \begin{picture}(20,40)(0,0)
    \put(10,20){\makebox(0,0){$=$}}
  \end{picture}
  \begin{picture}(40,40)(0,0)	
    \halfgrid
    \put(0,10){\line(1,0){30}}
    \put(0,30){\line(1,0){40}}
    \put(30,35){\line(1,0){10}}
    \put(10,0){\line(0,1){10}}
    \put(10,20){\line(0,1){20}}
    \put(30,0){\line(0,1){40}}
    \put(5,20){\line(0,1){10}}

    \put(7.5,23){\makebox(0,0)[b]{\scriptsize $a$}}
    \put(15,23){\makebox(0,0)[b]{\scriptsize $b$}}
    \put(25,13){\makebox(0,0)[b]{\scriptsize $c$}}
    \put(25,23){\makebox(0,0)[b]{\scriptsize $d$}}
  \end{picture}
  \begin{picture}(20,40)(0,0)
    \put(10,20){\makebox(0,0){$=$}}
  \end{picture}
  \begin{picture}(40,40)(0,0)	
    \halfgrid
    \put(0,10){\line(1,0){40}}
    \put(0,30){\line(1,0){40}}
    \put(10,0){\line(0,1){10}}
    \put(10,20){\line(0,1){20}}
    \put(30,0){\line(0,1){40}}
    \put(5,20){\line(0,1){10}}

    \put(7.5,23){\makebox(0,0)[b]{\scriptsize $a$}}
    \put(15,23){\makebox(0,0)[b]{\scriptsize $b$}}
    \put(25,13){\makebox(0,0)[b]{\scriptsize $c$}}
    \put(25,23){\makebox(0,0)[b]{\scriptsize $d$}}
  \end{picture}
\end{center}

\begin{center}
\setlength{\unitlength}{1.5pt}
\begin{picture}(40,40)(0,0)	
\end{picture}
\begin{picture}(20,50)(0,0)
  \put(10,20){\makebox(0,0){$=$}}
\end{picture}
\begin{picture}(40,40)(0,0)	
  \halfgrid
  \put(0,10){\line(1,0){40}}
  \put(0,30){\line(1,0){40}}
  \put(30,0){\line(0,1){10}}
  \put(30,20){\line(0,1){20}}
  \put(10,0){\line(0,1){40}}
  \put(5,20){\line(0,1){10}}

  \put(7.5,23){\makebox(0,0)[b]{\scriptsize $a$}}
  \put(15,23){\makebox(0,0)[b]{\scriptsize $b$}}
  \put(15,13){\makebox(0,0)[b]{\scriptsize $c$}}
  \put(25,23){\makebox(0,0)[b]{\scriptsize $d$}}
\end{picture}
\begin{picture}(20,40)(0,0)
  \put(10,20){\makebox(0,0){$=$}}
\end{picture}
\begin{picture}(40,40)(0,0)	
  \halfgrid
  \put(0,10){\line(1,0){40}}
  \put(10,30){\line(1,0){30}}
  \put(0,5){\line(1,0){10}}
  \put(30,0){\line(0,1){10}}
  \put(30,20){\line(0,1){20}}
  \put(10,0){\line(0,1){40}}
  \put(5,10){\line(0,1){10}}

  \put(7.5,13){\makebox(0,0)[b]{\scriptsize $a$}}
  \put(15,23){\makebox(0,0)[b]{\scriptsize $b$}}
  \put(15,13){\makebox(0,0)[b]{\scriptsize $c$}}
  \put(25,23){\makebox(0,0)[b]{\scriptsize $d$}}
\end{picture}
\begin{picture}(20,40)(0,0)
  \put(10,20){\makebox(0,0){$=$}}
\end{picture}
\begin{picture}(40,40)(0,0)	
  \halfgrid
  \put(0,10){\line(1,0){40}}
  \put(10,30){\line(1,0){30}}
  \put(0,5){\line(1,0){10}}
  \put(30,0){\line(0,1){40}}
  \put(10,0){\line(0,1){40}}

  \put(15,13){\makebox(0,0)[b]{\scriptsize $a$}}
  \put(15,23){\makebox(0,0)[b]{\scriptsize $b$}}
  \put(25,13){\makebox(0,0)[b]{\scriptsize $c$}}
  \put(25,23){\makebox(0,0)[b]{\scriptsize $d$}}
\end{picture}
\begin{picture}(20,40)(0,0)
  \put(10,20){\makebox(0,0){$=$}}
\end{picture}
\begin{picture}(40,40)(0,0)	
  \halfgrid
  \put(0,10){\line(1,0){40}}
  \put(0,30){\line(1,0){40}}
  \put(10,0){\line(0,1){40}}
  \put(30,0){\line(0,1){40}}
  \put(15,23){\makebox(0,0)[b]{\scriptsize $b$}}
  \put(25,23){\makebox(0,0)[b]{\scriptsize $d$}}
  \put(15,13){\makebox(0,0)[b]{\scriptsize $a$}}
  \put(25,13){\makebox(0,0)[b]{\scriptsize $c$}}
\end{picture}
\end{center}

\begin{center}
\setlength{\unitlength}{1.5pt}
\begin{picture}(40,40)(0,0)	
\end{picture}
\begin{picture}(20,50)(0,0)
  \put(10,20){\makebox(0,0){$=$}}
\end{picture}
\begin{picture}(40,40)(0,0)	
  \halfgrid
  \put(0,10){\line(1,0){40}}
  \put(0,30){\line(1,0){20}}
  \put(30,30){\line(1,0){10}}
  \put(20,5){\line(1,0){10}}
  \put(10,0){\line(0,1){40}}
  \put(30,0){\line(0,1){40}}
  \put(15,23){\makebox(0,0)[b]{\scriptsize $b$}}
  \put(25,13){\makebox(0,0)[b]{\scriptsize $d$}}
  \put(15,13){\makebox(0,0)[b]{\scriptsize $a$}}
  \put(25,6){\makebox(0,0)[b]{\scriptsize $c$}}
\end{picture}
\begin{picture}(20,40)(0,0)
  \put(10,20){\makebox(0,0){$=$}}
\end{picture}
\begin{picture}(40,40)(0,0)	
  \halfgrid
  \put(0,10){\line(1,0){40}}
  \put(0,30){\line(1,0){20}}
  \put(30,30){\line(1,0){10}}
  \put(20,5){\line(1,0){10}}
  \put(10,0){\line(0,1){10}}
  \put(10,20){\line(0,1){20}}
  \put(30,0){\line(0,1){40}}
  \put(35,10){\line(0,1){10}}
  \put(15,23){\makebox(0,0)[b]{\scriptsize $b$}}
  \put(32.5,13){\makebox(0,0)[b]{\scriptsize $d$}}
  \put(25,13){\makebox(0,0)[b]{\scriptsize $a$}}
  \put(25,6){\makebox(0,0)[b]{\scriptsize $c$}}
\end{picture}
\begin{picture}(20,40)(0,0)
  \put(10,20){\makebox(0,0){$=$}}
\end{picture}
\begin{picture}(40,40)(0,0)	
  \halfgrid
  \put(0,10){\line(1,0){40}}
  \put(0,30){\line(1,0){40}}
  \put(10,0){\line(0,1){10}}
  \put(10,20){\line(0,1){20}}
  \put(30,0){\line(0,1){40}}
  \put(35,10){\line(0,1){10}}
  \put(15,23){\makebox(0,0)[b]{\scriptsize $b$}}
  \put(32.5,13){\makebox(0,0)[b]{\scriptsize $d$}}
  \put(25,23){\makebox(0,0)[b]{\scriptsize $a$}}
  \put(25,13){\makebox(0,0)[b]{\scriptsize $c$}}
\end{picture}
\begin{picture}(20,40)(0,0)
  \put(10,20){\makebox(0,0){$=$}}
\end{picture}
\begin{picture}(40,40)(0,0)	
  \halfgrid
  \put(0,10){\line(1,0){40}}
  \put(0,30){\line(1,0){40}}
  \put(10,0){\line(0,1){40}}
  \put(30,0){\line(0,1){40}}
  \put(15,23){\makebox(0,0)[b]{\scriptsize $b$}}
  \put(25,23){\makebox(0,0)[b]{\scriptsize $a$}}
  \put(15,13){\makebox(0,0)[b]{\scriptsize $c$}}
  \put(25,13){\makebox(0,0)[b]{\scriptsize $d$}}
\end{picture}
\end{center}
(Note that the first eight moves effectuate a cyclic permutation of
the four middle elements, and that the last four moves permutes three
of the middle elements.)
\end{dem}

The remainder of this section explores a few consequences of the
commutativity phenomenon.

\begin{blanko}{Cancellative double semigroups.}
  A double semigroup is {\em cancellative} if $x * c = y * c$ implies 
  $x=y$, for multiplication with any $c$ from any of the four sides.
  The following corollary is immediate.
\end{blanko}

\begin{cor}\label{cancellative}
  A cancellative double semigroup is commutative.
  \qed
\end{cor}

\noindent
More generally it is sufficient that there exists one cancellable element
$c$ whose powers in both directions are also cancellable:
then place this $c$ in the empty boxes of the argument and cancel.

\begin{blanko}{Inverse double semigroups.}
  Recall that two elements $x$ and $y$ in a semigroup are said to be
  each other's inverses if $xyx=x$ and $yxy=y$, and that an {\em
  inverse semi\-group} is one where every element has a unique inverse.
  This can also be described as a universal algebraic structure: an
  associative multiplication, together with a unary `inverse'
  operation $x \mapsto x^{-1}$, satisfying $(xy)^{-1} = y^{-1}x^{-1}$,
  and $xx^{-1}yy^{-1} = yy^{-1}xx^{-1}$.  Inverse semi\-groups are
  important in many areas of mathematics and arise notably as
  semi\-groups of partial symmetries,
  cf.~Clifford-Preston~\cite{Clifford-Preston}.  See also
  Lawson~\cite{Lawson:InverseSemigroups}, who explains 
  the equivalence (due to Ehresmann and Schein) between inverse 
  semigroups and certain ordered groupoids.

  An {\em inverse double semigroup} is a double semigroup both of
  whose semigroup structures are inverse.  We shall see in a moment 
  that inverse double semigroups are commutative.
\end{blanko}

\begin{lemma}
  In an inverse double semigroup, the horizontal and vertical
  inverse operations commute.  
\end{lemma}

\begin{dem}
  Given an element $a\in S$, let $a^{-1}$ denote the horizontal 
  inverse, and let $\sigma(a)$ denote the vertical inverse.  We need 
  to show that $\sigma(a^{-1}) = \sigma(a)^{-1}$.
  
  Claim 1:
  \begin{center}
  \setlength{\unitlength}{1.5pt}
  \begin{picture}(44,30)(0,0)
    \put(0,0){\line(1,0){42}}
    \put(0,10){\line(1,0){42}}
    \put(0,20){\line(1,0){42}}
    \put(0,30){\line(1,0){42}}
    \put(0,0){\line(0,1){30}}
    \put(12,10){\line(0,1){10}}
    \put(30,10){\line(0,1){10}}
    \put(42,0){\line(0,1){30}}
    \put(21,4){\makebox(0,0)[b]{\scriptsize $a$}}
    \put(21,24){\makebox(0,0)[b]{\scriptsize $a$}}
    \put(6,13){\makebox(0,0)[b]{\scriptsize $\sigma(a)$}}
    \put(21,13){\makebox(0,0)[b]{\scriptsize $\sigma(a^{-1})$}}
    \put(37,13){\makebox(0,0)[b]{\scriptsize $\sigma(a)$}}
  \end{picture}
  \begin{picture}(20,30)(0,0)
    \put(10,15){\makebox(0,0){$=$}}
  \end{picture}
  \begin{picture}(10,30)(0,-10)	
    \put(0,10){\line(1,0){10}}
    \put(0,0){\line(1,0){10}}
    \put(0,0){\line(0,1){10}}
    \put(10,0){\line(0,1){10}}
    \put(5,4){\makebox(0,0)[b]{\scriptsize $a$}}
  \end{picture}
\end{center}

Claim 2:
  \begin{center}
  \setlength{\unitlength}{1.5pt}
  \begin{picture}(44,30)(0,0)
    \put(0,0){\line(1,0){42}}
    \put(0,10){\line(1,0){42}}
    \put(0,20){\line(1,0){42}}
    \put(0,30){\line(1,0){42}}
    \put(0,0){\line(0,1){30}}
    \put(12,0){\line(0,1){10}}
    \put(30,0){\line(0,1){10}}
    \put(42,0){\line(0,1){30}}
    \put(12,20){\line(0,1){10}}
    \put(30,20){\line(0,1){10}}
    \put(21,14){\makebox(0,0)[b]{\scriptsize $a$}}
    \put(6,3){\makebox(0,0)[b]{\scriptsize $\sigma(a)$}}
    \put(21,3){\makebox(0,0)[b]{\scriptsize $\sigma(a^{-1})$}}
    \put(37,3){\makebox(0,0)[b]{\scriptsize $\sigma(a)$}}
    \put(6,23){\makebox(0,0)[b]{\scriptsize $\sigma(a)$}}
    \put(21,23){\makebox(0,0)[b]{\scriptsize $\sigma(a^{-1})$}}
    \put(37,23){\makebox(0,0)[b]{\scriptsize $\sigma(a)$}}
  \end{picture}
  \begin{picture}(20,30)(0,0)
    \put(10,15){\makebox(0,0){$=$}}
  \end{picture}
  \begin{picture}(42,30)(0,-10)	
    \put(0,10){\line(1,0){42}}
    \put(0,0){\line(1,0){42}}
    \put(0,0){\line(0,1){10}}
    \put(12,0){\line(0,1){10}}
    \put(30,0){\line(0,1){10}}
    \put(42,0){\line(0,1){10}}
    \put(6,3){\makebox(0,0)[b]{\scriptsize $\sigma(a)$}}
    \put(21,3){\makebox(0,0)[b]{\scriptsize $\sigma(a^{-1})$}}
    \put(37,3){\makebox(0,0)[b]{\scriptsize $\sigma(a)$}}
  \end{picture}
\end{center}

Both claims follow by starting on the left-hand side by rewriting $a=aa^{-1}a$.
Then use the interchange law, and compute each column.

These two equations show that $\sigma(a)\sigma(a^{-1})\sigma(a)$ is 
the vertical inverse to $a$, in other words
$$
\sigma(a)\sigma(a^{-1})\sigma(a) = \sigma(a) .
$$

Now repeat the arguments with $a^{-1}$ and $a$ interchanged, to show 
also
$$
\sigma(a^{-1})\sigma(a)\sigma(a^{-1}) = \sigma(a^{-1}) .
$$
These two equations show that $\sigma(a^{-1})$ is the horizontal inverse
to $\sigma(a)$.  That is, $\sigma(a^{-1}) = \sigma(a)^{-1}$ as we 
wanted to show.
\end{dem}

In view of the lemma, we can adopt the following notation: given an
element $a$, let $a^{-1}$ denote the horizontal inverse, and let $A$
denote the vertical inverse.  Then by the lemma, the notation $A^{-1}$
is unambiguous.

\begin{prop}\label{inverse}
  Every inverse double semigroup is commutative.  
\end{prop}

\begin{dem}
  Let $A$ and $B$ be two elements in an inverse double semigroup.
  The proof consists in writing a big multiplication where $AB$ 
  appears in the middle, as inverse to the outer factors. Then
  commute $A$ and $B$ using Proposition~\ref{ab}.  By uniqueness of inverses
  we can then conclude that $AB=BA$.  In fact this argument is needed 
  four times: one for each way one element can be inverse to another.
  Here goes:
  
  Claim 1:
  \begin{center}
  \setlength{\unitlength}{1.5pt}
  \begin{picture}(44,30)(0,0)
    \put(0,0){\line(1,0){42}}
    \put(0,10){\line(1,0){42}}
    \put(0,20){\line(1,0){42}}
    \put(0,30){\line(1,0){42}}
    \put(0,0){\line(0,1){30}}
    \put(42,0){\line(0,1){30}}
    \put(21,3){\makebox(0,0)[b]{\scriptsize $ab$}}
    \put(21,23){\makebox(0,0)[b]{\scriptsize $ab$}}
    \put(21,13){\makebox(0,0)[b]{\scriptsize $ABA^{-1}B^{-1}AB$}}
  \end{picture}
  \begin{picture}(20,30)(0,0)
    \put(10,15){\makebox(0,0){$=$}}
  \end{picture}
  \begin{picture}(10,30)(0,-10)	
    \put(0,10){\line(1,0){10}}
    \put(0,0){\line(1,0){10}}
    \put(0,0){\line(0,1){10}}
    \put(10,0){\line(0,1){10}}
    \put(5,3){\makebox(0,0)[b]{\scriptsize $ab$}}
  \end{picture}
\end{center}

Claim 2:
  \begin{center}
  \setlength{\unitlength}{1.5pt}
  \begin{picture}(44,30)(0,0)
    \put(0,0){\line(1,0){42}}
    \put(0,10){\line(1,0){42}}
    \put(0,20){\line(1,0){42}}
    \put(0,30){\line(1,0){42}}
    \put(0,0){\line(0,1){30}}
    \put(42,0){\line(0,1){30}}
    \put(21,13){\makebox(0,0)[b]{\scriptsize $ab$}}
    \put(21,3){\makebox(0,0)[b]{\scriptsize $ABA^{-1}B^{-1}AB$}}
    \put(21,23){\makebox(0,0)[b]{\scriptsize $ABA^{-1}B^{-1}AB$}}
  \end{picture}
  \begin{picture}(20,30)(0,0)
    \put(10,15){\makebox(0,0){$=$}}
  \end{picture}
  \begin{picture}(42,30)(0,-10)	
    \put(0,10){\line(1,0){42}}
    \put(0,0){\line(1,0){42}}
    \put(0,0){\line(0,1){10}}
    \put(42,0){\line(0,1){10}}
    \put(21,3){\makebox(0,0)[b]{\scriptsize $ABA^{-1}B^{-1}AB$}}
  \end{picture}
\end{center}
 
Claim 1 and 2 together say that $ABA^{-1}B^{-1}AB$ is the vertical 
inverse to $ab$.  That is,
$$
AB \; A^{-1}B^{-1}\; AB = AB .
$$
Repeating the arguments with $( )^{-1}$ on every symbol we find 
also
$$
A^{-1}B^{-1}\; AB\; A^{-1}B^{-1} = A^{-1}B^{-1} .
$$
These two equation together say that $AB$ is the horizontal
inverse to $A^{-1}B^{-1}$.  But so is $BA$, and we are done,
provided we can prove the two claims.
  
Let us just prove Claim 1.  The 
second claim is analogous.
\begin{center}
\setlength{\unitlength}{1.5pt}
\begin{picture}(44,50)(0,-20)
  \put(0,0){\line(1,0){42}}
  \put(0,10){\line(1,0){42}}
  \put(0,20){\line(1,0){42}}
  \put(0,30){\line(1,0){42}}
  \put(0,0){\line(0,1){30}}
  \put(42,0){\line(0,1){30}}
  \put(21,3){\makebox(0,0)[b]{\scriptsize $ab$}}
  \put(21,23){\makebox(0,0)[b]{\scriptsize $ab$}}
  \put(21,13){\makebox(0,0)[b]{\scriptsize $ABA^{-1}B^{-1}AB$}}
\end{picture}
\begin{picture}(16,50)(0,0)
  \put(8,39){\makebox(0,0)[b]{\tiny (1)}}
  \put(8,35){\makebox(0,0){$=$}}
\end{picture}
\begin{picture}(44,50)(0,-20)
  \put(0,0){\line(1,0){42}}
  \put(0,10){\line(1,0){42}}
  \put(0,20){\line(1,0){42}}
  \put(0,30){\line(1,0){42}}
  \put(0,0){\line(0,1){30}}
  \put(42,0){\line(0,1){30}}
  \put(21,3){\makebox(0,0)[b]{\scriptsize $abb^{-1}a^{-1}ab$}}
  \put(21,23){\makebox(0,0)[b]{\scriptsize $abb^{-1}a^{-1}ab$}}
  \put(21,13){\makebox(0,0)[b]{\scriptsize $ABA^{-1}B^{-1}AB$}}
\end{picture}
\begin{picture}(16,50)(0,0)
  \put(8,39){\makebox(0,0)[b]{\tiny (2)}}
  \put(8,35){\makebox(0,0){$=$}}
\end{picture}
\begin{picture}(60,50)(0,0)
  \put(0,0){\line(1,0){60}}
  \put(0,10){\line(1,0){60}}
  \put(0,20){\line(1,0){60}}
  \put(0,30){\line(1,0){60}}
  \put(0,40){\line(1,0){60}}
  \put(0,50){\line(1,0){60}}
  \put(0,0){\line(0,1){50}}
  \put(10,0){\line(0,1){50}}
  \put(20,0){\line(0,1){50}}
  \put(30,0){\line(0,1){50}}
  \put(40,0){\line(0,1){50}}
  \put(50,0){\line(0,1){50}}
  \put(60,0){\line(0,1){50}}
  \put(5,3){\makebox(0,0)[b]{\tiny $a$}}
  \put(15,3){\makebox(0,0)[b]{\tiny $b$}}
  \put(25,3){\makebox(0,0)[b]{\tiny $b^{-1}$}}
  \put(35,3){\makebox(0,0)[b]{\tiny $a^{-1}$}}
  \put(45,3){\makebox(0,0)[b]{\tiny $a$}}
  \put(55,3){\makebox(0,0)[b]{\tiny $b$}}
  \put(5,13){\makebox(0,0)[b]{\tiny $A$}}
  \put(15,13){\makebox(0,0)[b]{\tiny $B$}}
  \put(25,13){\makebox(0,0)[b]{\tiny $B^{-1}$}}
  \put(35,13){\makebox(0,0)[b]{\tiny $A^{-1}$}}
  \put(45,13){\makebox(0,0)[b]{\tiny $A$}}
  \put(55,13){\makebox(0,0)[b]{\tiny $B$}}
  \put(5,23){\makebox(0,0)[b]{\scriptsize $a$}}
  \put(15,23){\makebox(0,0)[b]{\scriptsize $b$}}
  \put(25,23){\makebox(0,0)[b]{\scriptsize $b^{-1}$}}
  \put(35,23){\makebox(0,0)[b]{\scriptsize $a^{-1}$}}
  \put(45,23){\makebox(0,0)[b]{\scriptsize $a$}}
  \put(55,23){\makebox(0,0)[b]{\scriptsize $b$}}
  \put(5,33){\makebox(0,0)[b]{\scriptsize $A$}}
  \put(15,33){\makebox(0,0)[b]{\scriptsize $B$}}
  \put(25,33){\makebox(0,0)[b]{\scriptsize $A^{-1}$}}
  \put(35,33){\makebox(0,0)[b]{\scriptsize $B^{-1}$}}
  \put(45,33){\makebox(0,0)[b]{\scriptsize $A$}}
  \put(55,33){\makebox(0,0)[b]{\scriptsize $B$}}
  \put(5,43){\makebox(0,0)[b]{\scriptsize $a$}}
  \put(15,43){\makebox(0,0)[b]{\scriptsize $b$}}
  \put(25,43){\makebox(0,0)[b]{\scriptsize $b^{-1}$}}
  \put(35,43){\makebox(0,0)[b]{\scriptsize $a^{-1}$}}
  \put(45,43){\makebox(0,0)[b]{\scriptsize $a$}}
  \put(55,43){\makebox(0,0)[b]{\scriptsize $b$}}
  \linethickness{1.5pt}
  \put(19.5,40){\line(1,0){21}}
  \put(19.5,30){\line(1,0){21}}
  \put(20,29.5){\line(0,1){11}}
  \put(40,29.5){\line(0,1){11}}
  \put(9.5,50){\line(1,0){41}}
  \put(9.5,10){\line(1,0){41}}
  \put(10,9.5){\line(0,1){41}}
  \put(50,9.5){\line(0,1){41}}
\end{picture}
\begin{picture}(16,50)(0,0)
  \put(8,39){\makebox(0,0)[b]{\tiny (3)}}
  \put(8,35){\makebox(0,0){$=$}}
\end{picture}
\begin{picture}(60,50)(0,0)
  \put(0,0){\line(1,0){60}}
  \put(0,10){\line(1,0){60}}
  \put(0,20){\line(1,0){60}}
  \put(0,30){\line(1,0){60}}
  \put(0,40){\line(1,0){60}}
  \put(0,50){\line(1,0){60}}
  \put(0,0){\line(0,1){50}}
  \put(10,0){\line(0,1){50}}
  \put(20,0){\line(0,1){50}}
  \put(30,0){\line(0,1){50}}
  \put(40,0){\line(0,1){50}}
  \put(50,0){\line(0,1){50}}
  \put(60,0){\line(0,1){50}}
  \put(5,3){\makebox(0,0)[b]{\tiny $a$}}
  \put(15,3){\makebox(0,0)[b]{\tiny $b$}}
  \put(25,3){\makebox(0,0)[b]{\tiny $b^{-1}$}}
  \put(35,3){\makebox(0,0)[b]{\tiny $a^{-1}$}}
  \put(45,3){\makebox(0,0)[b]{\tiny $a$}}
  \put(55,3){\makebox(0,0)[b]{\tiny $b$}}
  \put(5,13){\makebox(0,0)[b]{\tiny $A$}}
  \put(15,13){\makebox(0,0)[b]{\tiny $B$}}
  \put(25,13){\makebox(0,0)[b]{\tiny $B^{-1}$}}
  \put(35,13){\makebox(0,0)[b]{\tiny $A^{-1}$}}
  \put(45,13){\makebox(0,0)[b]{\tiny $A$}}
  \put(55,13){\makebox(0,0)[b]{\tiny $B$}}
  \put(5,23){\makebox(0,0)[b]{\scriptsize $a$}}
  \put(15,23){\makebox(0,0)[b]{\scriptsize $b$}}
  \put(25,23){\makebox(0,0)[b]{\scriptsize $b^{-1}$}}
  \put(35,23){\makebox(0,0)[b]{\scriptsize $a^{-1}$}}
  \put(45,23){\makebox(0,0)[b]{\scriptsize $a$}}
  \put(55,23){\makebox(0,0)[b]{\scriptsize $b$}}
  \put(5,33){\makebox(0,0)[b]{\scriptsize $A$}}
  \put(15,33){\makebox(0,0)[b]{\scriptsize $B$}}
  \put(25,33){\makebox(0,0)[b]{\scriptsize $B^{-1}$}}
  \put(35,33){\makebox(0,0)[b]{\scriptsize $A^{-1}$}}
  \put(45,33){\makebox(0,0)[b]{\scriptsize $A$}}
  \put(55,33){\makebox(0,0)[b]{\scriptsize $B$}}
  \put(5,43){\makebox(0,0)[b]{\scriptsize $a$}}
  \put(15,43){\makebox(0,0)[b]{\scriptsize $b$}}
  \put(25,43){\makebox(0,0)[b]{\scriptsize $b^{-1}$}}
  \put(35,43){\makebox(0,0)[b]{\scriptsize $a^{-1}$}}
  \put(45,43){\makebox(0,0)[b]{\scriptsize $a$}}
  \put(55,43){\makebox(0,0)[b]{\scriptsize $b$}}
  \linethickness{1.5pt}
  \put(19.5,40){\line(1,0){21}}
  \put(19.5,30){\line(1,0){21}}
  \put(20,29.5){\line(0,1){11}}
  \put(40,29.5){\line(0,1){11}}
  \put(9.5,50){\line(1,0){41}}
  \put(9.5,10){\line(1,0){41}}
  \put(10,9.5){\line(0,1){41}}
  \put(50,9.5){\line(0,1){41}}
\end{picture}
\begin{picture}(16,50)(0,0)
  \put(8,39){\makebox(0,0)[b]{\tiny (4)}}
  \put(8,35){\makebox(0,0){$=$}}
\end{picture}
\begin{picture}(10,30)(0,-30)	
  \put(0,10){\line(1,0){10}}
  \put(0,0){\line(1,0){10}}
  \put(0,0){\line(0,1){10}}
  \put(10,0){\line(0,1){10}}
  \put(5,3){\makebox(0,0)[b]{\scriptsize $ab$}}
\end{picture}
\end{center}

Step (1) is to rewrite using $ab= ab(ab)^{-1}ab = abb^{-1}a^{-1}ab$.
In Step (2), the bottom row is expanded into three rows.  Step (3) is
the crucial commutation, justified by Proposition~\ref{ab}.  In
Step (4), each column is reduced to a single symbol, and the resulting
six-fold horizontal product is resolved.
\end{dem}

\section{Two-fold monoidal categories}
\label{Sec:cat}

\begin{blanko}{Double semigroups in $\Cat$.}
  As remarked, Proposition~\ref{ab} holds true for double semigroups
  in any monoidal category.  We shall now be concerned with double
  semi\-groups in $\Cat$; these are categories equipped with two
  strictly associative multiplication laws satisfying the strict
  interchange law.  Observe that Proposition~\ref{ab} holds also
  when the symbols $a$ and $b$ represent arrows.
  
  As in the previous section, in order to get a useful and generic
  statement out of Proposition~\ref{ab}, we need the presence
  of some cancellability.  One interesting case (the only one we
  consider) is when the two multiplication laws possess weak units,
  that is when we have a two-fold monoidal category, strictly
  associative and with strict interchange.
  Weak units are in particular weakly cancellable, which is what we
  need.  The weakness means that we do not get strict commutativity as
  conclusion, but we do get a symmetry, and in fact a degenerate one.
  We say that a symmetry $\sigma$ on a monoidal category is {\em
  degenerate} if for every object $X$ we have $\sigma_{X,X} =
  \id_{X\tensor X}$.  The result is this:

  \begin{prop}\label{deg-sym}
  Let $\CC$ be a two-fold monoidal category, strictly associative 
  and with strict interchange, but only weak units.
  Then $\CC$ is degenerate symmetric.
\end{prop}
The proof starts on the next page.

\bigskip

Note that degeneracy is stable under braided monoidal equivalence, and
in particular most braided monoidal categories are not equivalent to
degenerate symmetric ones.  Hence we get this corollary:

\begin{cor}
  Not every two-fold monoidal category is equivalent 
  to one with strict associativity and strict interchange.
  \qed
\end{cor}

From standard facts about braided categorical groups and homotopy 
$3$-types, as explained in the Introduction, we then get:

\begin{cor}\label{3-groupoids}
  Strictly associative one-object, one-arrow $3$-groupoids (but
  still with weak units) cannot realise all simply-connected
  homotopy $3$-types.
  \qed
\end{cor}
%
\end{blanko}


\begin{blanko}{Weak units.}
  A weak unit for a category with (strict) multiplication is an object
  $I$ equipped natural isomorphisms $\lambda_X: I \tensor X \isopil X$
  and $\rho_X : X \tensor I \isopil X$ satisfying $\id_X \tensor
  \lambda_Y = \rho_X \tensor \id_Y$.
  Alternatively~\cite{Kock:units1}, a weak unit can be characterised
  as an object $I$ equipped with a single isomorphism $\alpha: I \tensor I
  \isopil I$ and having the property that tensoring with it from
  either side is an equivalence of categories.  (The left and right
  constraints can be canonically constructed from $\alpha$,
  cf.~\cite{Kock:units1}.)  The property that tensoring with $I$ from
  either side is an equivalence is the crucial property for the
  present purposes --- we shall say that $I$ is {\em cancellable}.
\end{blanko}

\begin{blanko}{Two-fold strictly associative monoidal categories.}
  Let $(\CC, \tensor_h, \tensor_v)$ denote a strict double semigroup
  in $\Cat$, where each multiplication has a weak unit, denoted respectively
  $(I_h, \alpha_h)$ and $(I_v,\alpha_v)$. 
  Certain compatibilities could be required of the two structures,
  for example commutativity of the square \maltese\ in the next proof,
  but we shall not need any such further conditions.
%
%
%
\end{blanko}

\begin{lemma}
  The two units are isomorphic as objects in $\CC$: $I_h \simeq I_v$.
\end{lemma}
\begin{dem}
  In fact there are two different isomorphisms, depending on which 
  route we take in the diagram
  \begin{center}
  \setlength{\unitlength}{1.5pt}
  \begin{picture}(50,57)(0,0)
    
    \put(-35,5){
    \put(0,0){\line(0,1){10}}
    \put(10,0){\line(0,1){10}}
    \put(0,0){\line(1,0){10}}
    \put(0,10){\line(1,0){10}}
    \put(5,3){\makebox(0,0)[b]{\scriptsize $I_h$}}
    }
    
    \put(-12.5,10){\makebox(0,0)[b]{\scriptsize $\simeq$}}
    
    \put(0,5){
    \put(0,0){\line(0,1){10}}
    \put(10,0){\line(0,1){10}}
    \put(20,0){\line(0,1){10}}
    \put(0,0){\line(1,0){20}}
    \put(0,10){\line(1,0){20}}
    \put(5,3){\makebox(0,0)[b]{\scriptsize $I_h$}}
    \put(15,3){\makebox(0,0)[b]{\scriptsize $I_h$}}
    }
    
        \put(30,10){\makebox(0,0)[b]{\scriptsize $\simeq$}}

    \put(45,40){
    \put(0,0){\line(1,0){10}}
    \put(0,10){\line(1,0){10}}
    \put(0,20){\line(1,0){10}}
    \put(0,0){\line(0,1){20}}
    \put(10,0){\line(0,1){20}}
    \put(5,3){\makebox(0,0)[b]{\scriptsize $I_v$}}
    \put(5,13){\makebox(0,0)[b]{\scriptsize $I_v$}}
}

    \put(67.5,50){\makebox(0,0)[b]{\scriptsize $\simeq$}}

    \put(80,45){
    \put(0,0){\line(0,1){10}}
    \put(10,0){\line(0,1){10}}
    \put(0,0){\line(1,0){10}}
    \put(0,10){\line(1,0){10}}
    \put(5,3){\makebox(0,0)[b]{\scriptsize $I_v$}}
    }
    
    \put(32.5,50){\makebox(0,0)[b]{\scriptsize $\simeq$}}
    \put(10,27.5){\makebox(0,0)[b]{\scriptsize $\verticalSimeq$}}
    \put(50,30){\makebox(0,0)[b]{\scriptsize $\verticalSimeq$}}

    \put(30,30){\makebox(0,0){\maltese}}

    \put(0,40){
    
     \twobytwogrid
    \put(5,2){\makebox(0,0)[b]{\scriptsize $I_v$}}
    \put(15,2){\makebox(0,0)[b]{\scriptsize $I_h$}}
    \put(5,13){\makebox(0,0)[b]{\scriptsize $I_h$}}
    \put(15,13){\makebox(0,0)[b]{\scriptsize $I_v$}}
}

    \put(40,0){
    \twobytwogrid
    \put(5,2){\makebox(0,0)[b]{\scriptsize $I_h$}}
    \put(15,2){\makebox(0,0)[b]{\scriptsize $I_v$}}
    \put(5,13){\makebox(0,0)[b]{\scriptsize $I_v$}}
    \put(15,13){\makebox(0,0)[b]{\scriptsize $I_h$}}
    }
  \end{picture}
\end{center} 
  
\vspace{-2\bigskipamount}
\end{dem}

We shall not need any specific isomorphism --- just from its existence we 
conclude that $I_h$ is cancellable also vertically.   So 
tensoring with $I_h$ from any of the four sides is an 
equivalence of categories,  that's all we need to know.
We now set
$$
I \df I_h ,
$$
and forget about the vertical unit.

%
%
%

%


\begin{dem*}{Proof of Proposition~\ref{deg-sym}.}
  Step 1: the functor 
\begin{eqnarray*}
  F:\CC & \longrightarrow & \CC  \\
  X & \longmapsto & \begin{array}{ccc}
    I & I & I  \\
    I & X & I  \\
    I & I & I  \\
    I & I & I 
  \end{array}
\end{eqnarray*}
is an equivalence of categories.  This is clear since it is the composite of
tensoring with $I$ from each of the four sides.

Step 2: the functor $F$ has a strong multiplicative structure (with respect to
$\tensor\df\tensor_h$).  In other words, there are natural isomorphisms 
$FX\tensor FY \isopil F(X\tensor Y)$ satisfying the usual 
associativity condition \cite[XI.2]{MacLane:categories}.  Specifically,
$$
\begin{array}{cccccc}
  I & I & I & I & I & I  \\
  I & X & I & I & Y & I  \\
  I & I & I & I & I & I  \\
  I & I & I & I & I & I 
\end{array}
\isopil
\begin{array}{cccccc}
  I & I  & I  \\
  I & XY  & I  \\
  I & I  & I  \\
  I & I  & I 
\end{array}
$$
is defined row-wise as the composite of horizontal left and right constraints,
multiplying together horizontally the four middle columns.
By coherence for the horizontal unit, it does not matter how these
multiplications are effectuated.

Step 3: $F$ has the property that for any pair of objects (or arrows)
$X, Y \in \CC$, we have $FX\tensor FY = FY \tensor FX$.
This follows immediately from Proposition~\ref{ab}.

Now we define a symmetry $\sigma$ on $\tensor$ component-wise by
transporting these equalities back along the multiplicative
equivalence $F$.  Explicitly, by fully faithfulness and strong 
multiplicativity of $F$ we have bijections
$$
\CC( X\tensor Y,Y\tensor X) \simeq \CC(F(X\tensor Y),F(Y\tensor X)) \simeq \CC(FX\tensor FY, 
FY\tensor FX)
$$
and in this last hom set we have a distinguished element, namely the
identity arrow.  Define $\sigma_{X,Y}$ to be the arrow $X\tensor Y 
\to Y \tensor X$ corresponding
to the identity arrow under this bijection.  It is easy to see that
the $\sigma_{X,Y}$ are natural in $X$ and $Y$, and that they satisfy
the axioms for a symmetry
\cite[Def.~2.1]{Joyal-Street:braided-tensor}.  (Note that
compatibility with the unit is automatic from the axioms,
cf.~\cite[Prop.~2.1]{Joyal-Street:braided-tensor}.)   It is also
clear that $\sigma$ is degenerate.
\end{dem*}

\nocite{Kock:0507}
\nocite{Simpson:9810}
\nocite{Balteanu-et.al.:9808}
\nocite{Eckmann-Hilton}
\nocite{Joyal-Kock:traintracks}
\nocite{Kock:units1}
\nocite{MacLane:naturalAssociativity}
\nocite{MacLane:categories}
\nocite{Lawson:InverseSemigroups}

\footnotesize

\hyphenation{mathe-matisk}


\begin{thebibliography}{10}

\bibitem{Balteanu-et.al.:9808}
{\sc Cornel B{\u a}lteanu, Zbigniew Fiedorowicz, Roland Schw{\"a}nzl, {\rm and
  }Rainer Vogt}.
\newblock {\em Iterated monoidal categories}.
\newblock Adv. Math. {\bf 176} (2003), 277--349.
\newblock ArXiv:math.AT/9808082.


  \bibitem{Brown-Gilbert}
{\sc Ronald Brown {\rm and }Nicholas D.~Gilbert}.
\newblock {\em Algebraic models of {$3$}-types and automorphism structures for
  crossed modules}.
\newblock Proc. London Math. Soc. (3) {\bf 59} (1989), 51--73.

\bibitem{Cech:1932}
{\sc Eduard {\v{C}}ech}.
\newblock {\em H{\"o}herdimensionale Homotopiegruppen}.
\newblock In {\em Verhandlungen des internationalen Mathematiker-Kongresses,
  Z{\"u}rich 1932}, vol.~2, page 203. Walter Saxer, Z\"urich, 1932.
\newblock Reprint, Kraus, Nendeln, Liechtenstein, 1967.

\bibitem{Clifford-Preston}
{\sc Alfred H. Clifford {\rm and }Gordon B. Preston}.
\newblock {\em The algebraic theory of semigroups. {V}ol. {I+II}}.
\newblock No.~7 in Mathematical Surveys. American Mathematical Society,
  Providence, R.I., 1961;1967.

\bibitem{Eckmann-Hilton}
{\sc Beno Eckmann {\rm and }Peter J.~Hilton}.
\newblock {\em Group-like structures in general categories I: multiplication
  and comultipliplcation}.
\newblock Math. Ann. {\bf 145} (1962), 227--255.

\bibitem{Gordon-Power-Street}
{\sc Robert Gordon, A.~John Power, {\rm and }Ross Street}.
\newblock {\em Coherence for tricategories}.
\newblock Mem. Amer. Math. Soc. {\bf 117} (1995), vi+81pp.

\bibitem{Joyal-Kock:traintracks}
{\sc Andr{\'e} Joyal {\rm and }Joachim Kock}.
\newblock {\em Weak units and homotopy $3$-types}.
\newblock To appear in the StreetFest
  Proceedings (Contemp. Math., edited by M.~Batanin, A.~Davydov, M.~Johnson,
  S.~Lack, and A.~Neeman). ArXiv:math.CT/0602084, 

\bibitem{Joyal-Street:braided-tensor}
{\sc Andr{\'e} Joyal {\rm and }Ross Street}.
\newblock {\em Braided tensor categories}.
\newblock Adv. Math. {\bf 102} (1993), 20--78.
\newblock Formerly Macquarie Math. Reports No.~850067 (Dec.~1985) and
  No.~860081 (Nov.~1986).

\bibitem{Joyal-Tierney:homotopy-types}
{\sc Andr{\'e} Joyal {\rm and }Myles Tierney}.
\newblock {\em Algebraic homotopy types}.
\newblock Unpublished manuscript (1984).

\bibitem{Kock:0507}
{\sc Joachim Kock}.
\newblock {\em Weak identity arrows in higher categories}.
\newblock IMRP Internat. Math. Res. Papers {\bf 2006} (2006), 1--54.
\newblock ArXiv:math.CT/0507116.

\bibitem{Kock:units1}
{\sc Joachim Kock}.
\newblock {\em Elementary remarks on units in monoidal categories}.
\newblock Preprint, arXiv:math.CT/0507349.

\bibitem{Lawson:InverseSemigroups}
{\sc Mark~V. Lawson}.
\newblock {\em Inverse semigroups.  The theory of partial symmetries}.
\newblock World Scientific Publishing Co. Inc., River Edge, NJ, 1998.

\bibitem{Loday:catn}
{\sc Jean-Louis Loday}.
\newblock {\em Spaces with finitely many nontrivial homotopy groups}.
\newblock J. Pure Appl. Alg. {\bf 24} (1982), 179--202.

\bibitem{MacLane:naturalAssociativity}
{\sc Saunders {Mac~Lane}}.
\newblock {\em Natural associativity and commutativity}.
\newblock Rice Univ. Studies {\bf 49} (1963), 28--46.

\bibitem{MacLane:categories}
{\sc Saunders {Mac~Lane}}.
\newblock {\em Categories for the working mathematician, second edition}.
\newblock No.~5 in Graduate Texts in Mathematics. Springer-Verlag, New York,
  1998.

\bibitem{Paoli:0607}
{\sc Simona Paoli}.
\newblock {\em Semistrict models of connected 3-types and Tamsamani's weak
  3-groupoids}.
\newblock Preprint, arXiv:math.AT/0607330.

\bibitem{Simpson:9810}
{\sc Carlos Simpson}.
\newblock {\em Homotopy types of strict $3$-groupoids}.
\newblock Preprint, arXiv:math.CT/9810059.

\bibitem{Tamsamani:thesis}
{\sc Zouhair Tamsamani}.
\newblock {\em Sur des notions de $n$-cat\'egorie et $n$-groupo\"\i de non
  strictes via des ensembles multi-simpliciaux}.
\newblock $K$-Theory {\bf 16} (1999), 51--99.
\newblock (alg-geom/9512006 and alg-geom/9607010).

\end{thebibliography}


\noindent 
{\sc Departament de Matem\`atiques -- Universitat Aut\`onoma de Barcelona
-- 08193 Bellaterra (Barcelona) -- Spain}

\label{lastpage}

\end{document}